\documentclass[12pt,a4paper]{article}
\usepackage{amsmath, amssymb, theorem, latexsym}
\usepackage{color}

\newtheorem{theorem}{Theorem}[section]
\newtheorem{lemma}[theorem]{Lemma}
\newtheorem{proposition}[theorem]{Proposition}
\newtheorem{definition}[theorem]{Definition}

\newtheorem{remark}[theorem]{Remark}

\allowdisplaybreaks

\newcommand \clb{\color{blue}}
\newcommand\finbox{~\hfill$\Box$}%

\def\Om {{\Omega}}
\def\la {{\lambda}}

\def \Inte{{\rm Int\,}}

\newcommand {\pa}{\partial}

\numberwithin{equation}{section}
\begin{document}
%\null\vspace{4cm}
{\centering
\bfseries
{\Large Spectral minimal partitions for a thin strip on a cylinder or a thin annulus like domain with Neumann condition}\\

{\renewcommand\thefootnote{}%
\footnote{1991 Mathematics Subject Classification 35B05}}
\par
\mdseries
\scshape
\small
B. Helffer$^{1}$\\
T. Hoffmann-Ostenhof$^{2,3}$ \\\
\par
\upshape
Laboratoire de Math\'ematiques, Universit\'e  Paris-Sud 11$^{1}$\\
Institut f\"ur Theoretische Chemie, Universit\"at Wien$^2$\\
International Erwin Schr\"odinger Institute for Mathematical Physics$^3$\\
\date{}

}
\begin{abstract}
We analyze "Neumann"  spectral minimal partitions for a thin strip on a cylinder or   for the  thin annulus.
\end{abstract}
\section{Introduction}
In previous papers, sometimes in collaborations with other colleagues, we have analyzed spectral minimal partitions for some specific 
open subsets of $\mathbb R^2$ and for the  sphere $\mathbb S^2$. See \cite{HHOT} for some of the basic results and for more detailed
definitions.
In contrast with two-dimensional eigenvalue problems  for which 
a few examples exist where the eigenvalues and the eigenfunctions are explicitly known -- rectangles, the disk, sectors, the equilateral
triangle, $\mathbb S^2$ and the torus -- explicit non-nodal examples for minimal partitions are lacking. Up to now we only have been able 
to work out explicitly $\mathfrak L_3$ for the $2$-sphere. \cite{HHOT:2010}.  

Here we find other examples of non-nodal minimal partitions for problems for which the circle $\mathbb S^1$ is a deformation retract. Note that the Laplacian on the circle $\mathbb S^1_*$ (with perimeter 1)  can be interpreted as the Laplacian on  an interval $(0,1)$ with periodic boundary conditions. For this one dimensional problem we can work out 
the partition eigenvalues (see below for a definition) $\mathfrak L_k(\mathbb S_*^1)$ explicitly. We have $\mathfrak L_k=\pi^2k^2$. Observe that for odd $k\ge 3$ the 
$\mathfrak L_k$ are not eigenvalues, whereas for k even they are. 
The corresponding k-partitions are given by partitioning the circle into $k$ equal parts, hence  $D_1=(0,1/k), D_2=(1/k, 2/k), \dots , D_k=((k-1)/k, 1)$ identifying $0$ with $1$.  

We will consider a strip on a cylinder or  the annulus with suitable boundary conditions. All these domains are homotopic to 
$\mathbb S^1_*$.  For those domains we are going to investigate the corresponding minimal $3$-partitions. 

We recall some notation and definitions. Consider a k-partition $\mathcal D_k=(D_1,\dots, D_k)$ i.e. $k$ disjoint open subsets $D_i$ of 
some  $\Om$. Here $\Om$ can be a bounded domain in $\mathbb R^2$ or in a $2$-dimensional $C^\infty$ Riemannian manifold.

Consider first $-\Delta$ on $\Om$ where $\Delta$ can be the usual Laplacian or in the case of a manifold (with boundary or without boundary)  the corresponding 
Laplace-Beltrami operator. For the case with boundary we can impose  Dirichlet or  Neumann but we could also have mixed boundary conditions.

We associate with  $\mathcal D_k$ 
\begin{equation*}
 \Lambda(\mathcal D_k)=\sup_{1\le i\le k}  \la_1(D_i)
\end{equation*}
where   $\la_1(D_i)$ denotes
\begin{itemize}
\item either the lowest eigenvalue of the Dirichlet Laplacian in $D_i$
\item or the lowest eigenvalue of the Laplacian in $D_i$ where we put the Dirichlet boundary condition
on $\pa D_i\subset \Om$ and the Neumann boundary condition on  $\pa D_i\cap \pa \Om$.  \end{itemize}

It is probably worth to explain rigorously what we mean above by $\lambda_1(D_i)$ 
in the case of  measurable $D_i$'s.
\begin{definition}\label{eigenvalue}~\\
For any measurable $\omega \subset \Omega$, let $\lambda_1^D(\omega)$ (resp. $ \lambda_1^N(\omega)$) denotes the
first eigenvalue of the Dirichlet realization (resp. $\pa \Omega$-Neumann)  of the  operator in the following generalized sense. 
We define
$$
\lambda_1^{D\,or\, N} (\omega)=+\infty\;,
$$
if
$
\left\{ u\in W^1(\Omega)\,, u\equiv 0\; \text{a.e. on}\;
 \Omega\setminus \omega \right\}=\left\{0\right\}
$, 
$$
\lambda_1^D(\omega)=\inf \left\{\frac{\int_\Omega |\nabla
u(x)|^2 \, dx}{\int_\Omega|u(x)|^2\, dx}\;:\;
 u\in W^1_0(\Omega)\setminus\{0\}\,, u\equiv 0\; \text{a.e. on}\;
 \Omega\setminus \omega \right\}\;,
 $$
 $$
\lambda_1^N(\omega)=\inf \left\{\frac{\int_\Omega |\nabla
u(x)|^2 \, dx}{\int_\Omega |u(x)|^2\, dx}\;:\;
 u\in W^1(\Omega)\setminus\{0\}\,, u\equiv 0\; \text{a.e. on}\;
 \Omega\setminus \omega \right\}\;,
 $$
otherwise.\\
 We call groundstate any function $\phi$ achieving the above infimum. 
\end{definition}
Of course, if $\omega \subset\subset \Omega$, we have $\lambda_1^D(\omega)=\lambda_1^N(\omega)$.

The $k$-th partition-eigenvalue
$\mathfrak L_k(\Om)$ is then defined by 
\begin{equation}\label{frakLk}
 \mathfrak L_k(\Om)=\inf_{\mathcal D}\Lambda(\mathcal D)\,,
\end{equation}
where the infimum is considered\footnote{We refer to \cite{HHOT} for a more precise definition 
 of the considered class of $k$-partitions and the notion of regular representatives.} over the $k$-partitions.\\
Any $k$-partition $\mathcal D$ for which 
\begin{equation}\label{miniSP}
 \mathfrak L_k(\Om)=\Lambda (\mathcal D)
\end{equation}
is called spectral minimal $k$-partition, for short minimal $k$-partition. \\

If needed we will write $\mathfrak L_k^D(\Omega)$ or $\mathfrak L_k^N(\Omega)$ to indicate if we choose the Dirichlet condition or the $\pa\Omega$-Neumann condition in the above definitions .\\

Although not explicitly written in \cite{HHOT}, all the results obtained in the case of Dirichlet  are also true in the case of  Neumann. In particular, minimal partitions exist and have regular representatives.
 
One of the main results in \cite{HHOT} concerns the characterization of the case of equality in Courant's nodal Theorem. 
Consider an eigenvalue problem $-\Delta u_k=\la_k$ with suitable homogeneous boundary conditions (as previously defined)  and order 
the eigenvalues in increasing order $\la_1<\la_2\le \la _3 \le \dots \le  \la_k\dots$. If we assume that $u_k$ is  real, then Courant's nodal theorem says
that the number of its nodal domains   $\mu(u_k)$ satisfies  $\mu(u_k)\le k$. Note that Courant's nodal theorem holds in greater generality, in higher dimensions and with a potential.  Here a nodal domain is a component of $\Om\setminus N(u_k)$ where $\Om$ is the domain in $\Om$ or the manifold and $N(u_k)=\overline{\{x\in \Om\:|\: u_k(x)=0\}}.$ 
We call $u_k$ and $\la_k$ {\bf Courant sharp} if $\mu(u_k)=k$. In \cite{HHOT} we have also described some properties of minimal partitions. In many respects they 
are related to nodal domains. Nodal domains have many interesting properties. In particular in neighboring nodal domains the corresponding eigenfunction has different signs. Thereby two nodal domains $D_i, D_j$ are said to be neighbors if $\overline{\Inte{D_i\cup D_j}}$
is connected. We can associate with any (not necessarily nodal)  partition, say $\mathcal D_k=(D_1,\dots, D_k)$,
 a simple graph in the following way:  we associate to each $D_i$ a vertex and draw an edge between two vertices $i,j$ 
if the corresponding $D_i, D_j$ are neighbors. This amounts to 
say  that nodal graphs $\mathcal G(\mathcal D_k)$ are {\bf bipartite} graphs. 

The relation with Courant's nodal theorem is now the following, which is valid in the Dirichlet or Neumann case: 
\begin{theorem}[Dirichlet]\label{HHOTD}~\\
 If for a bounded domain $\Om$ with smooth boundary a minimal $k$-partition $\mathcal D$ with associated partition eigenvalue $\mathfrak L_k^D$  has a bipartite graph 
$\mathcal G(\mathcal D)$, then this minimal $k$-partition is produced by the nodal domains of an eigenfunction $u$ which is {\bf Courant sharp}
so that $-\Delta^{D} u =\lambda_k^{D} u $ in $\Om$ and $\la_k^{D} =\mathfrak L_k^{D}$. 
\end{theorem}
\begin{theorem}[$\pa\Omega$-Neumann]\label{HHOTN}~\\
 If for a bounded domain $\Om$ with smooth boundary a minimal $k$-partition $\mathcal D$ with associated partition eigenvalue $\mathfrak L_k^N$  has a bipartite graph 
$\mathcal G(\mathcal D)$, then this minimal $k$-partition is produced by the nodal domains of an eigenfunction $u$ which is {\bf Courant sharp}
so that $-\Delta^{N} u =\lambda_k^{N} u $ in $\Om$ and $\la_k^{N} =\mathfrak L_k^{N}$. 
\end{theorem}
Note that by the minimax principle $\la_k^D \le \mathfrak L_k^D$, resp. $\lambda_k^N\leq \mathfrak L_k^N$ and that, by Pleijel's Theorem \cite{Pleijel}, for each $\Om$ there is a $k(\Om)$ such that for each $k>k(\Om)$ any associated eigenfunction $u$ has strictly less than $k$ nodal domains. That implies that, for sufficiently high $k$,  the spectral  minimal   $k$-partitions are non-nodal.  \\
Note also that we will also meet mixed cases, when either Dirichlet or Neumann boundary conditions are assumed on different components of $\pa \Omega$.

\section{Neumann problem for a strip on the cylinder}
We start with the a strip $C(1,b)=C(b)  $ on a cylinder where 
\begin{equation}\label{C(b)}
C(b)=\mathbb S^1_* \times (0,b)\,.
\end{equation}
 If needed, we can represent the strip by a rectangle $R(1,b)=(0,1)\times (0,b)$
with identification of  $x=0$ and $x=1$. But the open sets of the partition are always  considered as open sets on the strip.

We consider Neumann boundary conditions at $y=0$ and $y=b$. The spectrum for the  Laplacian $\Delta^{N}$ with these  boundary conditions is  discrete and  is given by  
\begin{equation}\label{laC}
\sigma(-\Delta^{N})=\{ \pi^2 (4 m^2 + \frac{n^2}{b^2} )_{(m,n) \in \mathbb N^2}\}.
\end{equation}
Note that eigenvalues for $m\ge 1$ have at least multiplicity two. 
Identifying  $L^2(C(1,b))$ and  $L^2(R(1,b))$, a   corresponding orthonormal basis of eigenfunctions  is given by the functions on $R(1,b)$ 
$(x,y)\mapsto \cos (2\pi mx) \cos (\pi n \frac{y}{b})$ ($(m,n) \in \mathbb N^2$) and 
 $(x,y)\mapsto \sin (2\pi mx) \cos (\pi n \frac{y}{b})$ ($(m,n) \in \mathbb N^*\times \mathbb N$), where $\mathbb N^*=\mathbb N \setminus \{0\}$. We can now distinguish the following cases:
 \begin{enumerate} 
\item If $b<\frac 12$,
$$
\lambda_1^N =0\,, \,\lambda_2^N= \lambda_3^N= 4 \pi^2 < \lambda_4^N\,.
$$
\item If $\frac 12 < b < 1$, $$\lambda_1^N =0\,, \, \lambda_2^N= \frac{\pi^2}{b^2}\,,\, \lambda_3^N=\lambda_4^N = 4 \pi^2< \lambda_5^N \,.
$$
\item If $b=1$,  $$\lambda_1^N =0\,, \, \lambda_2^N= \pi^2\,, \, \lambda_3^N=\lambda_4^N=\lambda_5^N = 4 \pi^2 < \lambda_6^N\,.
$$
\item 
If $b>1$,
$$\clb 
\lambda_1^N =0\,,\, \lambda_2^N = \frac{ \pi^2}{b^2} \,,\, \lambda_3^N = \frac{4\pi^2}{b^2} < \lambda_4^N\,.
$$
\end{enumerate}
In particular, we see that $\lambda_3^N(C(1,b))$ is Courant sharp if and only if  $b\geq1$. Note also that, for $b\in (\frac 12,1]$, $\lambda_4^N(C(1,b)) $ cannot be Courant sharp,
 and that $\lambda_5^N(C(1,1))$ cannot be Courant sharp.

Before we state the main result for the strip on the cylinder,  we look also at its double covering $C(2,b)$, whose associated rectangle is given by $(0,2)\times (0,b)$. 
\begin{lemma}\label{C2}~\\
The Neumann eigenvalues for $C(2,b)$ are given, assuming that 
\begin{equation}\label{b<1/3}
b\le 1/3,
\end{equation}
by 
\begin{equation}\label{laC2}
 \la_1^N=0,\: \la_2^N=\la_3^N=\pi^2, \:\la_4^N=\la_5^N=4\pi^2, \:\la_6^N=\la_7^N=9\pi^2
\end{equation}
Note that $\la_6^N(C(2,b))$ is {\bf Courant sharp} if $b\le 1/3$. 
\end{lemma}
\begin{remark}~\\
Note that $\la_2^N=\la_3^N$ implies that $\mathfrak L_3^N>\la_3^N$ and that by Theorem \ref{HHOTN} the associated $\mathcal D_3$ is {\bf non-nodal}. \\
Note that for $b\geq 1$, we get by the same theorem that 
$\la_3^N(b)=\mathfrak L_3^N(b) = \frac{4 \pi^2}{b^2}$.  
\end{remark}

\begin{remark}\label{hitting rules}~\\
The Neumann boundary conditions imply that the zero's hit the boundary not as in the Dirichlet case. 
More precisely consider the cylinder $C(1,b)$ with associated rectangle $R(1,b)$ 
  and assume that the zero hits at  $(x_0,0)$. In polar coordinates $(r,\omega)$  centered at this point the zeroset looks locally like the zeroset of $r^m\cos m\omega,\: m=1,2,...$. This is in contrast with the Dirichlet case where we would have $r^m \sin m\omega$. The $r^m$ factor is just included to point out that 
eigenfunctions near zero's behave to leading order as harmonic homogeneous polynomials. 
\end{remark}
Here comes the main result for the
minimal $3$-partition for the strip on the cylinder.
\begin{theorem}\label{cylinder}~\\
 For 
\begin{equation}\label{b<!}
b\le b_0=\frac{1}{2\sqrt 5}\, ,
\end{equation}
we have 
\begin{equation}\label{9pi2}
 \mathfrak L_3^N(C(b))=9\pi^2. 
\end{equation}
The associated  minimal $3$-partition $\mathcal D_3(b)=(D_1,D_2,D_3)$ is  up to rotation represented by 
\begin{equation}\label{D3C}
D_\ell =((\ell -1)/3,\ell /3)\times (0,b),
\end{equation}  
in $R(1,b)$. 
\end{theorem}

Before giving the proof it might be appropriate to consider the case   $b_0<b\le 1$. 
 In this direction, we have :
\begin{proposition}\label{b>2/3}~\\
 For $b\in [2/3,1)$ the spectral minimal $3$-partition $\mathcal D_3(b)$ is not the one given by \eqref{D3C} and 
$\mathfrak L_3^N(C(b))<9\pi^2$ for $2/3<b<1$. 
\end{proposition}
 {\bf Proof of the proposition.}\\
It is immediate that the eigenfunction associated with $m=0$ and $n=2$  $(x,y) \mapsto \cos (2\pi \frac{y}{b})$ has three nodal domains with energy $\frac{4\pi^2}{b^2}$ which is less than $9\pi^2$.
\finbox

{\bf Proof of the theorem.}\\
Note that by definition  of $\mathfrak L_3^N$, we have in any case
\begin{equation}\label{L3<}
 \mathfrak L_3^N(C(b))\le 9\pi^2.
\end{equation}

We first sketch the main ideas for the proof. There are two arguments which will be crucial for the proof:\\
{\bf (1)} Take any candidate for a minimal 3-partition, 
$\mathcal D_3=(D_1,D_2,D_3)$. If we can show that $\Lambda(\mathcal D_3)>9\pi^2$ then this $\mathcal D_3$ cannot be a minimal partition due to the definition of $\mathfrak L_3^N(C(b))$.\\
{\bf (2)} Assume $b\le 1/3$. A 3-partition $\mathcal D_3$ is said to have {\bf property B} if  it becomes on the  double covering $C(2,b)$ a 6-partition.

Assume that there is {\bf minimal} 3-partition $\mathcal D_3$ with {\bf property B}.  Then $\Lambda(\mathcal D_3)$, the energy of this
partition is larger than or equal to  $\la_6^N(C(2,b))$. To see this just note that by Lemma \ref{C2},
$\la_6^N(C(2,b))$ is Courant sharp. The corresponding minimal $6$-partition $\mathcal D_6=(D_1,\dots, D_6)$ is given by 
\begin{equation}\label{calD6}
D_\ell=(\,(\ell-1)/3\,,\, \ell/3)\times (0,b),\:\:\: \ell=1,2,\dots, 6.
\end{equation}
Furthermore $\mathcal D_6$ for $C(2,b)$ is just the lifted 3-partition of $C(1, b)$ given in Theorem \ref{cylinder}.\\ 
{\bf Hence it suffices to show that the candidates for minimal partitions have the property B}. 

We also observe that the only candidates for minimal 3-partitions, assuming $b<1$, are non-nodal and further that, if we have 
a candidate for $\mathcal D_3=(D_1,D_2,D_3)$ for a minimal partition, each $D_i$ is nice, that means
\begin{equation}\label{nice}
\Inte(\overline D_i)=D_i.
\end{equation}
If not we could lower the energy by removing an arc inside $\overline D_i$ without reducing the number of the $D_i$. 
Here we neglect sets of capacity $0$. 

The assumption on $D_i$ implies by monotonicity that $$\la_1^D(D_i)> \lambda_1^{D}(C(1,b) ) = \pi^2/b^2.$$ Hence if $\pi^2/b^2>9\pi^2$ the associated partition must already lead to a $\Lambda_3(\mathcal D_3)>9\pi^2$, so {\bf(1)} applies. \\
{\bf End of the proof}\\
We proceed by showing that any minimal $3$-partition has  property {\bf B}. To show this it suffices that 
in any minimal $3$-partition $\mathcal D_3=(D_1,D_2,D_3)$ all the $D_i$ are 0-homotopic. This implies that lifting this partition
to the double covering yields a $6$-partition and the argument {\bf (2)} applies. 
We assume for contradiction that $D_3$ is not 0-homotopic, hence contains a path of index 1. We first observe that $D_1$ and $D_2$ must be neighbours (if not the partition would be nodal). Then let us introduce 
$D_{12} = \Inte(\bar D_1\cup \bar D_2)$. Because $D_3$ contains a path of index $1$,  $\overline{D_{12} }$ cannot touch one component  of the boundary of the cylinder and we have $\lambda_2^N(D_{12}) =\mathfrak L_3^N$. 
 By domain monotonicity  (this is not the  Dirichlet monotonicity result but the proof can be done either by reflection or by a density argument), 
the second eigenvalue $\lambda_2^N (D_{12} ) = \lambda_1^N(D_1)$ must be  be larger  than  the second eigenvalue of the Dirichlet-Neumann problem of the cylinder. But we have,  with $\la_i^{ND}$ denoting the eigenvalues with Neumann and Dirichlet boundary conditions  on the two components of the boundary of the strip, 
 \begin{equation}\label{laND}
 \la_1^{ND}=\frac{\pi^2}{4b^2}, \:\la_2^{ND}=\pi^2\min \Big(\frac{1}{b^2},\: \frac{1}{4b^2}+4\Big).
\end{equation}
Hence  
\begin{equation}
 \la_2^N(D_{1,2})>\la_2^{ND}\,,
\end{equation}
and we get a contradiction if $\la_2^{ND}\geq 9\pi^2$. We just have to work out the condition on $b$ such that 
\begin{equation}
 \min \big(\frac{1}{b^2},\:\frac{1}{4 b^2}+4\big)\geq 9.
\end{equation}
This is achieved for $b\leq (2\sqrt 5)^{-1}$ as claimed in  \eqref{b<!} in Theorem \ref{cylinder}. 

 \begin{remark}~\\
 We recall that, although there is a natural candidate (which is nodal on the double covering),  the minimal $3$-partition problem with Dirichlet  conditions for the annulus (also in the case of a thin annulus)  or the disk is still open.
 \end{remark}
 
 \begin{remark}~\\
 In view of the considerations above and of Proposition \ref{b>2/3} the question arises whether for some $b<1$ the corresponding 
minimal partition has the property that it has one or two points in its zero set where 3 arcs meet, hence having locally
a {\bf Y}-structure as discussed for instance in  \cite{HHOT:2010}.  As in the case of the rectangle considered in \cite{BHHO}, we can observe that in the case $b=1$, the eigenfunction $\cos 2\pi x
 - \cos 2 \pi y$ has a nodal set described in $R(1,1)$ by the two diagonals of the square
  and that it determines indeed a nodal $3$-partition. The guess is then that when $1-\epsilon< b<1$ (with $\epsilon >0$ small enough) this  nodal $3$-partition will  be deformed into  a non nodal minimal $3$-partition keeping the symmetry $(x,y) \mapsto (x, 1-y)$. We expect  two critical points from which three arcs start.
\end{remark}
 
\section{Extension to  minimal $k$-partitions of $C(b)$}
One can also consider  for $\Omega = C(b)$ minimal $k$-partitions with $k$ odd ($k\geq 3$) and Neumann condition and assume 
\begin{equation}\label{1k}
b < \frac 1k\,.
\end{equation}
The theorem of the previous section  can be extended to the case $k>3$. \\
First one observes that, if  the closure of one open set of the minimal $k$-partition contains a line joining the two components
 of the boundary, then one can go to the double covering and obtain a $(2k)$-partition. If  \eqref{1k} is satisfied, $\lambda_k^N(C(0,2b))$ is Courant sharp, and  get  as in the previous section  that the energy of this partition is necessarily higher than $k^2\pi^2$.
 
So there is no $D_i$ whose boundary has nonempty intersection
with both parts of the boundary of the strip. Hence there exists  one component of $\pa \Omega$  and  at least $\frac{k+1}{2}$ $D_i$ of the $k$-partition such that their boundaries  $\pa D_i$'s  do not intersect with this component. We immediately deduce that if  $b\leq \frac 1k$:
 \begin{equation}\label{imp}
 \min \left(\mathfrak L_\frac{k+1}{2} ^{DN} (\Omega), \mathfrak L_\frac{k+1}{2} ^{ND}(\Omega)\right) \geq k^2 \pi^2\,,
 \end{equation}
 we have $\mathfrak L_k^N (\Omega) = k^2 \pi^2$.\\
 Here $\mathfrak  L_\ell^{DN}$ corresponds to the $\ell$-th spectral partition eigenvalue for the strip with
Dirichlet boundary condition on $y=0$ and Neumann boundary condition for $y=b$ and $\mathfrak  L_\ell^{ND}$ is defined by exchange of the boundary conditions on the two boundaries.  In our special case, due to the symmetry with respect to $y=\frac b2$, we have actually $\mathfrak  L_\ell^{DN}(C(b)) =\mathfrak  L_\ell^{ND}(C(b))$.\\
 Having in mind that\footnote{Note that we have equality for $\ell$ even and $b <\frac{1}{\ell}$.} $\lambda_\ell^{DN} \leq  \mathfrak  L_\ell^{DN}$, 
 \eqref{imp} is a consequence of
 $$
 \la_\frac{k+1}{2}^{DN} \geq k^2 \pi^2.
 $$
 If $b <  \frac 1k$, we get in the case when $k= 4 p + 3$ ($p\in \mathbb N$) the additional condition 
 $$
 \frac {1}{4b^2} + 4 (p+1)^2 \geq  (4p+3)^2\,.
 $$
  Similarly, we get in the case when $k=4p+1$ ($p\in \mathbb N^*$)
  $$
 \frac {1}{4b^2} + 4 p^2  \geq (4p+1)^2\,.
 $$ 
 We have consequently proven:
 \begin{theorem}~\\
 If 
 \begin{itemize}
 \item
 $k=4p+3$ ($p\in \mathbb N$) and $b \leq 1\big / \sqrt{(3k+1)(k-1)}$,
 \item[]  or
 \item 
  $k=4p+1$ ($p\in \mathbb N^*$) and $b \leq 1/ \sqrt{(3k-1)(k+1)}$,
 \end{itemize} then $$ \mathfrak L_k(C(b)) = k^2\pi^2\,,$$ and a minimal $k$-partition is given by  $D_\ell=((\ell-1)/k,\ell/k))\times (0,b)\, $,\,  for $\ell =1,\dots, k$.
 \end{theorem}
 
\section{Generalization to other thin domains}
The previous proof is more general than it seems at the first look. At the price of less explicit results 
 we have a similar result for an annulus like domain $\Omega$. 
 We mention first  the case k=3 where the conditions read

 \begin{itemize}
 \item The eigenfunction associated with $\lambda_6^N(\Omega^R)$ is Courant sharp and antisymmetric with respect to the deck transformation from $\Omega^R$ onto $\Omega$,  
 \item 
 \begin{equation}\label{condthin}
 \lambda_6^N(\Omega^R) \leq \inf (\lambda_2^{DN}(\Omega), \lambda_2^{ND}(\Omega))\,.
 \end{equation}
 \end{itemize}
 
 Here $\Omega^R$ is the double covering of $\Omega$ and (DN) (respectively (ND))  corresponds
 to the Dirichlet-Neumann problem (Dirichlet inside, Neumann outside), respectively  (Dirichlet outside, Neumann inside). 
In the case of the annulus, these conditions can be made more explicit.

Here is a typical result which can be expected. For $b >0$ and two regular functions $h_1(\theta)$ and $h_2(\theta)$ on the circle such that $h_1 < h_2$, we consider an annulus like domain around the unit circle defined in polar coordinates by
$$
A(b)= \{(x,y) \,:\,  1+b  h_1(\theta) < r < 1+ b h_2(\theta)\} \,.
$$
It is clear from \cite{KZ} together with Poincar\'e's inequality that there exists $b_0 >0$ such that, if $0 <b\leq b_0$, condition \eqref{condthin} 
is satisfied. \\
One should also show the condition for Courant sharpness, which is true for the sixth eigenvalue of the lifted  Laplacian on the double covering of the annulus and  should be
 also  true  for our more general situation but for which we have no references, (see however \cite{FK} for thin curved tubes and \cite{L}).

\begin{remark}~\\
Although not explicit, condition \eqref{condthin} can be analyzed by perturbative method. This is indeed a purely spectral question. There is a huge literature  concerning thin domains, see for example
 \cite{FS, KZ} (and references therein).
\end{remark}

\begin{remark}~\\
Similar considerations lead also to extensions to higher $k$ odd for the thin annulus with Neumann boundary conditions. 
\end{remark}

{\bf Acknowledgements}~\\
The second author acknowledges helpful discussions with Frank Morgan during the 2010 Dido conference in Carthage.

\scshape 
B. Helffer: Laboratoire de Math\'ematique, B\^at. 425, Universit\'e  Paris Sud 11, 91 405 Orsay (France). \\  e-mail address : Bernard.Helffer@math.u-psud.fr.

\scshape
T. Hoffmann-Ostenhof: Institut f\"ur Theoretische Chemie, Universit\"at
Wien, W\"ahringer Strasse 17, A-1090 Wien, Austria and International Erwin
Schr\"odinger Institute for Mathematical Physics, Boltzmanngasse 9, A-1090
Wien, Austria.\\
 e-mail address : thoffman@esi.ac.at.

\end{document}